\documentclass[draft]{amsart}
\usepackage{amssymb} 
\newtheorem{thm}{Theorem}[section]
\newtheorem{cor}[thm]{Corollary}
\newtheorem{lemma}[thm]{Lemma}
\newtheorem{pro}[thm]{Proposition}

\newtheorem{exa}[thm]{Example}
\newtheorem{qe}[thm]{Question}
\newtheorem{Note}[thm]{Notification}

\newtheorem{rem}[thm]{Remark}
\numberwithin{equation}{section}

\begin{document}

\title[$L^T_p-$ functions on  locally compact groups]
{$L^T_p-$ functions on  locally compact groups}

\author[F. Abtahi H. G. Amini and A. Rejali]
{F. Abtahi,$^1$$^{*}$, H. G. Amini,$^2$ and A. Rejali,$^3$}





\date{Received: xxxxxx; Revised: yyyyyy; Accepted: zzzzzz.
\newline \indent $^{*}$Corresponding author}

\begin{abstract}
Let $G$ be a locally compact group and $1\leq p<\infty$. Based on
some important earlier works, in this paper the concept of
$L_p^T-$function is introduced. Then the structure of the space
$L^{T}_p(G)$, which is consisting of all $L_p^T-$functions, is
investigated. As an important result, $L^{T}_p(G)^+$ is completely
characterized, for the class of amenable locally compact groups.
Furthermore, it is verified some applications of Plancherel
Theorem on $L^{T}_2(G)$, in the class of abelian locally compact
groups.
\end{abstract} \maketitle

\section{\bf Introduction and basic definitions}

Let $G$ be a locally compact group and $\lambda$ be a fixed left
Haar measure on $G$. For $1\leq p<\infty$, the Lebesgue space
$L^p(G)$ with respect to $\lambda$ is as defined in \cite{hr};
i.e. the equivalence classes of $\lambda-$measurable functions on
$G$ with
$$
\|f\|_p=\left(\int_G|f(x)|^pd\lambda(x)\right)^{1/p}.
$$
We denote this space by $\ell^p(G)$, in the case where $G$ is
discrete. Not that for each $1\leq p<\infty$, $L^p(G)\subseteq
L^1(G)+L^{\infty}(G)$. In fact for each $f\in L^p(G)$, let
$$
A=\{x\in G: |f(x)|\leq 1\}.
$$
Then $f=f\chi_A+f\chi_{G\setminus A}$, where $f\chi_A\in
L^{\infty}(G)$ and $f\chi_{G\setminus A}\in L^1(G)$.

For $\lambda-$measurable functions $f$ and $g$ on $G$, the
convolution multiplication
$$
f*g(x)=\int_Gf(y)g(y^{-1}x)d\lambda(y)
$$
is defined at each point $x\in G$ for which this makes sense; i.e.
the function
$$
y\mapsto f(y)g(y^{-1}x)
$$
is $\lambda-$integrable. Then $f*g$ is said to exist if $f*g(x)$
exists for almost all $x\in G$. Note that the support of any given
function $f$ in $L^p(G)$ is contained in a $\sigma-$compact subset
of $G$. Thus if $f*g$ exists, for some $f,g\in L^p(G)$, then
\cite[Theorem 13.9]{hr} implies that $f*g$ is in fact a measurable
function on $G$.

According to \cite{g3}, \cite{m1} and \cite{m2}, for $1<p<\infty$,
a function $f\in L^p(G)$ is said to be $p-$tempered if
$$
\sup\{\|g*f\|_p,\ g\in C_{00}(G),\ \|g\|_p\leq 1\}<\infty,
$$
where $C_{00}(G)$ denotes the set of all continuous complex-valued
functions on $G$ with compact support. The set of all $p-$tempered
functions has been denoted by $L^t_p(G)$. Also the function
$\|.\|_p^t$ on $L^t_p(G)$ defined by
$$
\|f\|_p^t=\sup\{\|g*f\|_p,\ g\in C_{00}(G),\ \|g\|_p\leq 1\}
$$
is a norm, under which $L^t_p(G)$ is a normed vector space. It has
been proved in the first part of \cite[Theorem 1]{m1} that
\begin{eqnarray*}
\|f\|_p^t&=&\sup\{\|g*f\|_p,\ g\in C_{00}(G),\ \|g\|_p\leq 1\}\\
&=&\sup\{\|g*f\|_p,\ g\in L^1(G)\cap L^p(G),\ \|g\|_p\leq 1\}.
\end{eqnarray*}
The concept of $p-$tempered functions, was introduced and studied
in the decade of the 1970s; see \cite{g3}, \cite{m1}, \cite{m2}
and \cite{m3}. Mainly, Feichtinger \cite{f2} introduced the
general concept of tempered elements in Banach function spaces.
Indeed, let
\begin{enumerate}
\item[(i)] $(B,\|.\|)$ be an essential two-sided Banach module
over a Banach algebra $(A,\|.\|_A)$ with bounded two-sided
approximate identity; \item[(ii)] $A$ and $B$ are both
continuously embedded into a topological vector space and $A\cap
B$ is dense in $A$ as well as $B$; \item[(iii)] On $A\cap B$ the
multiplication inherited from $A$ and the two module operations
coincide.
\end{enumerate}
Then an element $c\in B$ is called tempered if
\begin{equation}\label{e21}
\|c\|^t=\sup\{\|dc\|: d\in A\cap B;\ \|d\|\leq 1\}<\infty.
\end{equation}
The space of all such elements will be denoted by $B^t$. For each
$c\in B^t$, the map $d\mapsto dc$ $(d\in A\cap B)$ can be extended
to a bounded linear operator on $B$ and it will be denoted by
$S_c$. Moreover, by \cite[Theorem 1.1]{f2}, $(B^t,\|.\|^t)$ is a
normed algebra with the multiplication defined by $a\circ
b=S_b(a)$, where $a,b\in B^t$. Furthermore, $(B,\||.\||^t)$ is a
Banach algebra, where $\||.\||^t=\|.\|+\|.\|^t$.\\

The concept of $p-$tempered in the $L^p-$spaces, can be taken from
this general definition of tempered, in the case where $G$ is a
unimodular locally compact group. Indeed, for $1\leq p<\infty$,
all the $L^p-$spaces can be continuously embedded in
$L^1_{loc}(G)$, the space of all locally integrable functions on
$G$. Now by takeing $B=L^p(G)$ and $A=L^1(G)$, the definition is
obtained. In this case, by \cite[Theorem 1.1]{f2},
$(L^t_p(G),\|.\|_p^t)$ is a normed algebra under convolution
product. Moreover $L^t_p(G)$ is a Banach algebra under convolution
and the following norm
$$
\||.\||_p^t=\|.\|_p+\|.\|^t_p.
$$

In \cite[section 2]{f2}, Feichtinger considered the solid Banach
function spaces $B\subseteq L^1_{loc}(G)$, which satisfy the
conditions (i), (ii) and (iii), whenever $A$ is a suitable
Beurling algebra $L^1(G,\omega)$, introduced in \cite{re}.
Moreover $B$ satisfies an additional condition as the following:
\begin{equation}\label{e14}
B*B\subseteq L^1_{loc}(G),
\end{equation}
in the sense that for all $f,g\in B$, $f*g$ exists as a function
on $G$, and belongs to $L^1_{loc}(G)$. Then in \cite[Theorem
2.5]{f2}, he showed that $f\in B$ is tempered if and only if $g*f$
is in $B$, for all $g\in B$. Moreover he proved that the extended
multiplication on $B^t$ is given by the ordinary convolution.

By \cite{f2}, if $G$ is unimodular and $1<p\leq 2$, then the
inclusion (\ref{e14}) is satisfied for $B=L^p(G)$. Also some other
sufficient conditions for $L^p(G)$ to satisfy \eqref{e14}, has
been provided in \cite[Page 55]{f2}. Thus in this case, $f\in
L^t_p(G)$ if and only if $g*f$ exists and belongs to $L^p(G)$, for
all $g\in L^p(G)$ and
$$
\|f\|_p^t=\sup\{\|g*f\|_p,\ g\in L^p(G),\ \|g\|_p\leq 1\}<\infty.
$$
\begin{Note}\label{n}\rm
We found further important remarks concerning $p-$tempered
functions in \cite{f2} and \cite{m3}, which are appropriate to
mention here. We draw the attention of the readers to some
noticeable points in the earlier papers. The concept of
$p-$tempered in the papers \cite{g3}, \cite{m1}, \cite{m2} and
\cite{m3} has been presented in two forms. One of them is that,
which we described before; i.e. $f\in L^p(G)$ is $p-$tempered if
\begin{equation}\label{e3}
\sup\{\|g*f\|_p,\ g\in C_{00}(G),\ \|g\|_p\leq 1\}<\infty.
\end{equation}
Another form is that $f\in L^p(G)$ is $p-$tempered if for all
$g\in L^p(G)$, $g*f$ exists and belongs to $L^p(G)$ and,
\begin{equation}\label{e20}
\|f\|_p^t=\sup\{\|g*f\|_p,\ g\in L^p(G),\ \|g\|_p\leq 1\}<\infty.
\end{equation}
In fact in \cite[Theorem 1]{m1}, it has been proved that the two
definitions are equivalent. It is clear that \eqref{e20} implies
\eqref{e3}, immediately. But as it is mentioned in \cite[Page
56]{f2} and \cite{m3}, there are some unclear points in the proof
of the converse of this statement. For compact groups,
\cite[Theorem 1]{m1} is obviously valid. Moreover \cite[Theorem
1]{m1} is valid in the case where $G$ is a [SIN]-group and $1<
p\leq 2$; see \cite{m3}. In fact the definition of tempered
elements, presented in \cite{f2}, is not identical with the
original one, presented in \cite{m2}. But as we mentioned before,
they coincide with each other by \cite[Theorem 2.5]{f2}. For other
important results on this topic, we just cite to \cite{f2},
\cite{g3}, \cite{m1}, \cite{m2} and \cite{m3}. See also \cite{e1}
and \cite{o} for the more general case of vector valued
$p-$tempered functions and weighted case, respectively.
\end{Note}

On the basis of the above discussions, we define the concept of
$L^T_p-$function, for each $1<p<\infty$, according to \eqref{e20}.
But in order to avoid any confusion with the original definitions
of $p-$tempered and $L_p^t$, we shall denote it with a different
symbol. In fact we call a function $f\in L^p(G)$, $L^T_p-$function
if $g*f$ exists, for all $g\in L^p(G)$ and also
\begin{equation}\label{e13}
\|f\|_p^T=\sup\{\|g*f\|_p,\ g\in L^p(G),\ \|g\|_p\leq 1\}<\infty.
\end{equation}
It is remarkable to mention that by \cite[Theorem 1.3]{r}, $f\in
L^{T}_p(G)$ if and only if for all $g\in L^p(G)$, $g*f$ exists and
belongs to $L^p(G)$. In fact only this statement implies the
condition \eqref{e13}. The set of all $L^T_p-$functions will be
denoted by $L^{T}_p(G)$. Moreover $\|.\|_p^T$ is a norm under
which $L^{T}_p(G)$ is a normed vector space. It is easily verified
that $L^{T}_p(G)\subseteq L^{t}_p(G)$ and for each $f\in
L^T_p(G)$, we have $\|f\|^T_p=\|f\|^t_p$. More importantly, in the
case where $G$ is unimodular we have
$$
L^1(G)\cap L^p(G)\subseteq L^{T}_p(G)\subseteq L^{t}_p(G).
$$
Another main point is that $L^{T}_p(G)$ is not in general a normed
algebra under convolution product. In fact Fubini's Theorem can
not be used and so the associativity rule is not satisfied. But by
the explanations given above, if $G$ is unimodular and $1<p\leq
2$, then
$$L^{t}_p(G)=L^{T}_p(G).$$ In this situation,
$(L^{T}_p(G),\|.\|_p^T)$ is a normed algebra under convolution
product.\\

Some other explanations about the construction of $L^{T}_p(G)$ is
appreciated. Let $G$ be a discrete group. Then the duality
characterization of $L^p$ implies that the positive functions of
$L^{T}_p(G)$, belong to $\ell^q(G)$. It should be noted that
$L^{T}_p(G)$ is not in general a subspace of $L^q(G)$. For
example, suppose that $G$ is a unimodular non discrete locally
compact group, and $1<p\leq 2$. Then
$$L^1(G)\cap L^p(G)\subseteq L^{T}_p(G)=L^{t}_p(G).$$ It follows
that $L^{T}_p(G)$ can not be contained in $L^q(G)$.

By the explanations given in Note \ref{n}, it seems that in
general, the concept of $L^T_p-$function is different from
$p-$tempered. However at the present time, we are not able to
provide examples to show the difference. We also could not find
any examples for this purpose, in the earlier publications. But
our results in the present work are self-contained and are not
coincided with the previous results.

The aim of the present paper is to investigate the structure of
the space $L^{T}_p(G)$. We first study some elementary properties
of this space. We also study $L^T_p(G)$ under another norm,
defined as the following
$$
\||f\||_p^T=\|f\|_p+\|f\|_p^T\;\;\;\;\;\;\;\;\;\;\;\;\;\;\;\;(f\in
L^T_p(G)).
$$
Furthermore, we compare three norms $\|.\|_p$, $\|.\|_p^T$ and
$\||.\||_p^T$ on $L^{T}_p(G)$, to show the different treatment of
$L^{T}_p(G)$, with respect to each norm. Then we follow the
purpose which this research has mainly raised. The main motivation
for this aim, stem from \cite[Lemma 2.3]{f2} which provides a
characterization of $L^t_p(G)^+$. Indeed, it has been proved that
every function $f$ in $L^t_p(G)^+$, the set of positive-valued
functions of $L^t_p(G)$, belongs to $L^1(G)\cap L^p(G)$, whenever
$G$ is a unimodular amenable locally compact group. Here, we
generalize this result and in fact by a completely self-contained
proof, we show that if $G$ is an amenable locally compact group,
then the set $L^T_p(G)^+$ is exactly the positive functions of
$L^1(G,\omega_q)\cap L^p(G)$, where $q$ is the exponential
conjugate of $p$, $\Delta$ is the modular function on $G$ and
$\omega_q=\Delta^{-1/q}$. Also $L^1(G,\omega_q)$ is the set of all
complex valued functions $f$ with $f\omega_q\in L^1(G)$. In the
last section, we present some applications of the Plancherel
Theorem on $L^{T}_2(G)$, whenever $G$ is an abelian locally
compact group. Using these achievements, we develop the Plancherel
Theorem for $L^{T}_2(G)\cap L^{\infty}(G)$ and show that the
restricted of the Plancherel mapping is an isometric isomorphism
onto $L^{T}_2(\widehat{G})\cap L^{\infty}(\widehat{G})$, under
some suitable norms.

\section{\bf Primary Results of $L^{T}_p(G)$}

In the present section, we study some general properties of
$L^{t}_p(G)$ and $L^{T}_p(G)$. This verification helps us to
obtain the main results of this paper. We commence with the
following lemma. Recall that $\delta_x$ is Dirac measure at $x$,
for each $x\in G$.
\begin{lemma}\label{l}
Let $G$ be a locally compact group and $1<p<\infty$. Then for each
$f\in L_p^{T}(G)$ (resp. $L_p^{t}(G)$) and $x\in G$,
$f*\delta_x,\delta_x*f\in L_p^{T}(G)$ (resp. $L_p^{t}(G)$).
Moreover
$\|f*\delta_x\|_p^T=\|\delta_x*f\|_p^T=\Delta(x)^{-1/q}\|f\|_p^T$
(resp.
$\|f*\delta_x\|_p^t=\|\delta_x*f\|_p^t=\Delta(x)^{-1/q}\|f\|_p^t$).
\end{lemma}

\begin{proof}
Let $f\in L_p^{T}(G)$ and $x\in G$. Then \cite[Theorem 20.13]{hr}
implies that for each $g\in L^p(G)$
\begin{eqnarray*}
\|g*(f*\delta_x)\|_p&=&\|(g*f)*\delta_x\|_p\\
&\leq&\|g*f\|_p\int_G
\Delta(y)^{-1/q}d\delta_x(y)\\
&\leq&\|g\|_p\|f\|_p^T\Delta(x)^{-1/q}.
\end{eqnarray*}
Thus $f*\delta_x\in L^{T}_p(G)$ and
$$
\|f*\delta_x\|_p^T\leq\|f\|_p^T\Delta(x)^{-1/q}.
$$
Consequently
$$
\|f\|_p^T=\|f*\delta_x*\delta_{x^{-1}}\|_p^T\leq\|f*\delta_{x}\|_p^T\Delta(x)^{1/q}
$$
and so $\|f*\delta_x\|_p^T=\Delta(x)^{-1/q}\|f\|_p^T$. Similarly,
one can show that $\delta_x*f\in L^{T}_p(G)$ and
$\|\delta_x*f\|_p^T=\Delta(x)^{-1/q}\|f\|_p^T$. Using a slight
modification of the above proof, we obtain the same result for
$f\in L_p^{t}(G)$.
\end{proof}
If $G$ is unimodular, then for $f\in L^1(G)\cap L^p(G)$ and each
$g\in L^p(G)$,
$$
\|g*f\|_p\leq\|g\|_p\|f\|_1
$$
and so
$$\|f\|_p^T\leq\|f\|_1.$$ In fact $L^1(G)\cap L^p(G)\subseteq
L_p^{T}(G)$. The converse can be established by adding an extra
condition.
\begin{pro}\label{p1}
Let $G$ be a locally compact group and $1<p<\infty$. Then the
following assertions are equivalent.
\begin{enumerate}
\item[(i)] $L^1(G)\cap L^p(G)\subseteq L_p^{T}(G)$ and
$\|f\|_p^T\leq \|f\|_1$, for each $f\in L^1(G)\cap L^p(G)$;
\item[(ii)] $G$ is unimodular.
\end{enumerate}
\end{pro}
\begin{proof}
We just prove that $(i)\Rightarrow (ii)$. Let $f\in L^1(G)\cap
L^p(G)$ be nonzero. By the hypothesis and Lemma \ref{l},
$$
\Delta(x)^{-1/q}\|f\|_p^T=\|\delta_x*f\|_p^T\leq
\|\delta_x*f\|_1=\|f\|_1.
$$
Thus
$$
\Delta(x)^{-1/q}\leq \frac{\|f\|_1}{\|f\|_p^T}
$$
and so $\Delta$ is bounded below from zero. It concludes that $G$
is unimodular.
\end{proof}
The space $L^{T}_p(G)$ is just $L^p(G)$, whenever $G$ is compact.
In fact $L^p(G)$ is closed under convolution if and only if $G$ is
compact; see \cite{r}. This is completely related to an ancient
conjecture, called '$L^p-$conjecture' which was proved finally by
Saeki \cite{s}, in the general case. The conjecture asserts that
$f*g$ exists and belongs to $L^p(G)$, for all $f,g\in L^p(G)$, if
and only if $G$ is compact. In the next result, some other
equivalent conditions to the equality of $L^p(G)=L_p^{T}(G)$ is
provided. Before, we recall the following theorem, which is
obtained from \cite[Theorems 1.1 and 2.5]{f2}.
\begin{thm}\label{t3}
Let $G$ be a unimodular locally compact group and $1<p\leq 2$.
Then $L^{T}_p(G)=L^{t}_p(G)$ and it is a Banach algebra under
$\||.\||_p^T$ and convolution product.
\end{thm}

\begin{thm}\label{t1}
Let $G$ be a locally compact group and $1<p<\infty$. Consider the
following assertions.
\begin{enumerate}
\item[(i)] $L^p(G)=L_p^{T}(G)$, as two sets. \item[(ii)] The norms
$\|.\|_p$ and $\||.\||_p^T$ are equivalent on $L^T_p(G)$.
\item[(iii)] There exists a positive constant $K$ such that
$\|f\|_p^T\leq K\|f\|_p$, for each $f\in L_p^{T}(G)$. \item [(iv)]
$G$ is compact.
\end{enumerate}
Then $(i)\Leftrightarrow (iv)$, $(ii)\Leftrightarrow (iii)$ and
$(iv)\Rightarrow (iii)$. If $G$ is unimodular and $1<p\leq 2$ then
all the assertions are equivalent.
\end{thm}

\begin{proof}
$(i)\Rightarrow (iv)$. Let $L^p(G)=L_p^{T}(G)$. Then $L^p(G)$ is
an algebra under convolution. So \cite[Theorem 1.3]{r} and
\cite{s} imply that $G$ is compact.

$(iv)\Rightarrow (i)$. It is clear from \cite{s}.

$(iii)\Rightarrow (ii)$. Suppose that there exists a positive
constant $K$ such that $\|f\|_p^T\leq K\|f\|_p$, for each $f\in
L_p^{T}(G)$. Thus
$$
\|f\|_p\leq\||f\||_p^T\leq (K+1)\|f\|_p
$$
and so two norms are equivalent on $L^{T}_p(G)$.

$(ii)\Rightarrow (iii)$. It is clear.

$(iv)\Rightarrow (iii)$. If $G$ is compact, then again by
\cite{s}, $L^p(G)$ is a Banach algebra and so for all $f,g\in
L^p(G)$
$$
\|g*f\|_p\leq\|g\|_p\|f\|_p.
$$
It follows that $\|f\|_p^T\leq\|f\|_p$ and so $(iii)$ obviously
holds.

Now let $G$ be unimodular and $1<p\leq 2$. We show that all the
assertions are equivalent. It is sufficient to prove
$(ii)\Rightarrow (i)$. Let the norms $\|.\|_p$ and $\||.\||_p^T$
be equivalent on $L^{T}_p(G)$. Then $L^{T}_p(G)$ is a Banach
algebra under the norm $\||.\||_p^T$ by Theorem \ref{t3}. It
follows that $L^{T}_p(G)$ is complete under $\|.\|_p$, as well.
Now the result follows by the density $L_p^{T}(G)$ in $L_p(G)$.
\end{proof}
In the next theorem, some equivalent conditions to the
completeness of $L^{T}_p(G)$ under the norm $\|.\|_p^T$ are
provided.
\begin{thm}\label{t2}
Let $G$ be a locally compact group and $1<p<\infty$. Consider the
following assertions.
\begin{enumerate}
\item[(i)] $L_p^{T}(G)$ is Banach algebra under $\|.\|_p^T$ and
convolution product. \item[(ii)] There exists a positive constant
$K$ such that $\||f\||_p^T\leq K\|f\|_p^T$, for each $f\in
L_p^{T}(G)$. \item[(iii)] The norms $\|.\|_p^T$ and $\||.\||^T_p$
are equivalent. \item[(iv)] There exists a positive constant $K$
such that $\|f\|_p\leq K\|f\|_p^T$, for each $f\in L_p^{T}(G)$.
\item[(v)] $G$ is discrete.
\end{enumerate}
Then $(v)\Rightarrow (iv)$, $(iv)\Rightarrow (iii)$ and
$(iii)\Rightarrow (ii)$. If $G$ is unimodular and $1<p\leq 2$ then
all the assertions are equivalent.
\end{thm}

\begin{proof}
$(v)\Rightarrow (iv)$. Let $G$ be discrete. By normalizing Haar
measure on $G$, we have
$$
\|f\|_p=\|\delta_e*f\|_p\leq\|f\|_p^T,
$$
and so the mentioned inequality in $(iv)$ is obtained by choosing
$K=1$.

$(iv)\Rightarrow (iii)$. If there exists a positive constant $K$
such that $\|f\|_p\leq K\|f\|_p^T$, for each $f\in L_p^{T}(G)$,
then
$$
\|f\|_p^T\leq\||f\||_p^T\leq (K+1)\|f\|_p^T
$$
and so $(iii)$ is obtained.

$(iii)\Rightarrow (ii)$. It is clear.

Now let $G$ be unimodular and $1<p\leq 2$. We show that all the
assertions are equivalent. It suffices to prove $(ii)\Rightarrow
(i)$ and $(i)\Rightarrow (v)$.

$(ii)\Rightarrow (i)$. Suppose that there exists a positive
constant $K$ such that $\||f\||_p^T\leq K\|f\|_p^T$, for each
$f\in L_p^{T}(G)$. Then two norm $\|.\|_p^T$ and $\||.\||_p^T$ are
equivalent on $L^{T}_p(G)$ and so the result is concluded by
Theorem \ref{t3}.

$(i)\Rightarrow (v)$. Let $L_p^{T}(G)$ be a Banach algebra under
$\|.\|_p^T$ and consider the identity map
$$
\iota: (L_p^{T}(G),\||.\||_p^T)\rightarrow (L_p^{T}(G),\|.\|_p^T).
$$
Since $\iota$ is continuous, open mapping theorem implies that
both norms $\|.\|_p^T$ and $\||.\||^T_p$ are equivalent. Thus
there exists a positive constant $K$ such that $\||f\||_p^T\leq
K\|f\|_p^T$, for each $f\in L_p^{T}(G)$. Thus
\begin{equation}\label{e10}
\|f\|_p\leq (K-1)\|f\|_p^T.
\end{equation}
Since $G$ is unimodular, then \cite[20.14]{hr} implies that for
each $g\in C_{00}(G)$, we have
\begin{equation}\label{e11}
\|g\|_p^T\leq\|g\|_1.
\end{equation}
For the completeness, we show that $L^1(G)\subseteq L^p(G)$. Let
$g\in L^1(G)$ and $(g_n)$ be a sequence in $C_{00}(G)$ such that
$\lim_{n\rightarrow\infty}\|g_n-g\|_1=0$. Inclusions (\ref{e10})
and (\ref{e11}) imply that $(g_n)$ is a cauchy sequence in
$L^p(G)$ and so it converges to an element $h\in L^p(G)$. Hence
$g=h$ almost every where on $G$. It follows that $g\in L^p(G)$ and
so $L^1(G)\subseteq L^p(G)$. Therefore $G$ is discrete by
\cite[Theorem 1]{z}.
\end{proof}

The following result is automatically fulfilled from Theorems
\ref{t1} and \ref{t2}.

\begin{cor}\label{c}
Let $G$ be a unimodular locally compact group and $1<p\leq 2$. The
norms $\|.\|_p$ and $\|.\|_p^T$ are equivalent on $L^{T}_p(G)$ if
and only if $G$ is finite.
\end{cor}

\begin{rem}\label{r}\rm
Let $G$ be a locally compact group and $1\leq p<\infty$.
\begin{enumerate}
\item[(i)] It is known that $L^1(G)$ is a Banach algebra with a
two-sided contractive approximate identity
$(e_{\alpha})_{\alpha\in\Lambda}$. It follows that
$L^1(G)=L^{t}_1(G)=L^{T}_1(G)$; indeed, for each $f\in L^1(G)$, we
have
$$
\|f\|_1=\lim_{\alpha\in\Lambda}\|e_{\alpha}*f\|_1\leq\|f\|_1^T=\sup_{\|g\|_1\leq
1}\|g*f\|_1\leq\|f\|_1.
$$
Thus $\|f\|_1=\|f\|_1^t=\|f\|_1^T$ and consequently $L^1(G)$,
$L^{t}_1(G)$ and $L^{T}_1(G)$ are isometrically isomorphic.
\item[(ii)] If $G$ is compact and $1<p<\infty$, then
$L_p^{t}(G)=L_p^{T}(G)=L^p(G)$ and so it possesses an approximate
identity which is bounded in $L^1(G)$. \item[(iii)] If
$1<p<\infty$, then $L^1(G)\subseteq L_p^{T}(G)$ if and only if $G$
is discrete. To see, first suppose that $G$ is discrete and
$f\in\ell^1(G)$. For each $g\in\ell^p(G)$ we have
$$
\|g*f\|_p\leq\|g\|_p\|f\|_1.
$$
Since $\ell^1(G)\subseteq\ell^p(G)$, it follows that $f\in
L^{T}_p(G)$ and $\|f\|_p^T\leq\|f\|_1$. The converse is obvious
from the fact that $L^1(G)\subseteq L^p(G)$ is equivalent to the
discreteness of $G$ \cite[Theorem 1]{z}. \item[(iv)] Let $G$ be a
unimodular locally compact group and $1<p\leq 2$. Then
$L_p^{t}(G)=L_p^T(G)$ and it is a left and right $L^1(G)-$module
by \cite[Theorem 1.1]{f2}.
\end{enumerate}
\end{rem}

It is known that $L^1(G)$ possesses a left (right) identity if and
only if $G$ is discrete. Since $L^p(G)$ is not generally an
algebra, when $1<p<\infty$, we define the concept of left quasi
identity for $L^p(G)$, as an element $e\in L^p(G)$ such that for
each $f\in L^p(G)$, $e*f\in L^p(G)$ and $e*f=f$. Right quasi
identity element is defined similarly. An element $e\in L^p(G)$ is
called quasi identity for $L^p(G)$ if it is both left and right
quasi identity.

\begin{pro}\label{p2}
Let $G$ be a locally compact group and $1<p<\infty$. Then $L_p(G)$
has a left (right) quasi identity if and only if $G$ is discrete.
\end{pro}
\begin{proof}
Let $e\in L_p(G)$ be a left quasi identity element for $L^p(G)$
and $U$ be a symmetric and relatively compact neighborhood of the
identity element of $G$. Hence there is a positive constant $K$
such that for each $x\in U$, $\Delta(x)^{-1/q}\leq K$. Suppose on
the contrary that $G$ is not discrete. Thus for each $n\in
{\mathbb{N}}$, there is a neighborhood $U_n\subseteq U$ such that
$\lambda(U_n)<1/n$. Thus for each $n\in {\mathbb{N}}$,
$e*\chi_{U_n}=\chi_{U_n}$ and so
\begin{eqnarray*}
\lambda(U_n)^{1/p}&=&\|\chi_{U_n}\|_p\\
&=&\|e*\chi_{U_n}\|_p\\
&\leq&\|e\|_p\int_G\Delta(x)^{-1/q}\chi_{U_n}(x)dx\\
&\leq&K\|e\|_p\lambda(U_n).
\end{eqnarray*}
It follows that for each $n\in {\mathbb{N}}$
$$
\|e\|_p\geq\frac{\lambda(U_n)^{1/p-1}}{K}\geq \frac{n^{1-1/p}}{K},
$$
which is a contradiction. Therefore $G$ is discrete. A similar
argument can be used with the assumption of the existence of a
right quasi identity. The converse is trivial, since the Dirac
measure at the identity element of $G$ is a quasi identity for
$\ell^p(G)$.
\end{proof}

\begin{pro}\label{p3}
Let $G$ be a unimodular locally compact group and $1<p\leq 2$.
Then $L_p^{T}(G)$ has an identity if and only if $G$ is discrete.
\end{pro}

\begin{proof}
If $G$ is discrete then $\delta_e$ is an identity for $L_p^{T}(G)$
and
$$
\|\delta_e\|_p^T=\sup\{\|f*\delta_e\|_p:\ f\in L^p(G),\
\|f\|_p\leq 1\}=1.
$$
For the converse suppose that $e\in L_p^{T}(G)$ is an identity for
$L_p^{T}(G)$. Thus for each $f\in L_p^{T}(G)$
$$
\|f\|_p=\|e*f\|_p\leq\|e\|_p\|f\|_p^T
$$
and Theorem \ref{t2} implies that $G$ is discrete.
\end{proof}

\section{\bf Main Results}

Feichtinger \cite{f2} and Gilbert \cite{g1} investigated the
structure of $L^{t}_p(G)$, whenever $G$ is unimodular and
amenable. In this section, we characterize the structure of
$L^{T}_p(G)^+$, in the case where $G$ is an amenable locally
compact group. We commence with following proposition.

\begin{pro}\label{p4}
Let $G$ be a locally compact group, $1<p<\infty$ and $f\in
L_p^{T}(G)$. Then the real valued functions $Im(f)$, $Re(f)$
belong to $L_p^{T}(G)$.
\end{pro}
\begin{proof}
First let $g\in L^p(G)$ with $Im(g)=0$. Then $Re(g*f)=g*Re(f)$ and
$Im(g*f)=g*Im(f)$. It follows that $g*Re(f)$ and $g*Im(f)$ belong
to $L^p(G)$. Also
\begin{equation}\label{e4}
\|g*Re(f)\|_p\leq\|g*f\|_p\leq\|g\|_p\|f\|_p^T
\end{equation}
and
\begin{equation}\label{e5}
\|g*Im(f)\|_p\leq\|g*f\|_p\leq\|g\|_p\|f\|_p^T.
\end{equation}
Now let $g\in L^p(G)$ be arbitrary. By the inequalities (\ref{e4})
and (\ref{e5}) we clearly obtain
$$
\|g*Re(f)\|\leq
2\|g\|_p\|f\|_p^T\;\;\;\;\;\;\;and\;\;\;\;\;\;\;\;\;\|g*Im(f)\|_p\leq
2\|g\|_p\|f\|_p^T.
$$
It follows that $Re(f),Im(f)\in L_p^{T}(G)$. Moreover
$$
\|Re(f)\|_p^T\leq 2\|f\|_p^T \;\;\;\;\;and\;\;\;\;\;
\|Im(f)\|_p^T\leq 2\|f\|_p^T.
$$
\end{proof}
\begin{rem}\label{r2}\rm
With regard to Proposition \ref{p4}, it may be expected that the
positive parts of $Re(f)$ and $Im(f)$ belong to $L^{T}_p(G)$, as
well. Contrary to the impression, in general the answer is vague
and inconclusive to us. However this expectation is realized,
whenever $f\in L^1(G,\omega_q)\cap L_p^{T}(G)$. It is in fact a
direct consequence of the following lemma.
\end{rem}

\begin{lemma}\label{l2}
Let $G$ be a locally compact group and $1<p<\infty$. Then
$$L^1(G,\omega_q)\cap L^p(G)\subseteq L_p^{T}(G).$$
\end{lemma}

\begin{proof}
Let $f\in L^1(G,\omega_q)\cap L^p(G)$. Then \cite[20.14]{hr}
implies that for each $g\in L^p(G)$,
$$
\|g*f\|_p\leq\|g\|_p\int_G\Delta(x)^{-1/q}|f(x)|dx=\|g\|_p\int_G|f(x)|\omega_q(x)dx.
$$
It follows that $f\in L^{T}_p(G)$ and
$$
\|f\|_p^T\leq\int_G|f(x)|\omega_q(x)dx.
$$
Therefore $L^1(G,\omega_q)\cap L^p(G)\subseteq L_p^{T}(G)$.
\end{proof}

We state here the main theorem of this section.

\begin{thm}\label{t4}
Let $G$ be an amenable locally compact group and $1<p<\infty$.
Then $L^1(G,\omega_q)\cap L^p(G)^+=L_p^{T}(G)^+$.
\end{thm}

\begin{proof}
Let $f\in L_p^{T}(G)^+$. Set
$\widetilde{f}=\Delta^{-1/p}\check{f}$, where
$\check{f}(x)=f(x^{-1})$, for all $x\in G$. It follows that
$\widetilde{f}\in L^p(G)$ and for all $k\in L^p(G)$ and
$(\widetilde{k}*f)^{\widetilde{}}=\widetilde{f}*k$. Hence
$\widetilde{f}*k\in L^p(G)$. Also we have
\begin{equation}\label{e6}
\|\widetilde{f}*k\|_p=\|\widetilde{k}*f\|_p\leq\|k\|_p\|f\|_p^T.
\end{equation}
The amenability of $G$ implies that for an arbitrary compact
subset $C$ of $G$ and $0<\varepsilon<1$, there exists a compact
subset $K$ of $G$ such that for each $x\in C$
\begin{equation}\label{e7}
\frac{|xK\cap K|}{|K|}>1-\varepsilon,
\end{equation}
see \cite{p}. Take
$$
g=\frac{\chi_K}{|K|^{1/p}}\;\;\;\;\;\;\;and\;\;\;\;\;h=\frac{\chi_K}{|K|^{1/q}}.
$$
The inequality (\ref{e7}) and Fubini's Theorem imply that
\begin{eqnarray*}
\langle\widetilde{f}*g,h\rangle&=&\int_G\widetilde{f}(y)\frac{|yK\cap
K|}{|K|}dy\\
&\geq&\int_C\widetilde{f}(y)\frac{|yK\cap
K|}{|K|}dy\\
&\geq&(1-\varepsilon)\int_C\widetilde{f}(y)dy.
\end{eqnarray*}
This together with (\ref{e6}) yield that
$$
(1-\varepsilon)\int_C\widetilde{f}(y)dy\leq\langle\widetilde{f}*g,h\rangle\leq\|g\|_p\|f\|_p^T\|h\|_q=\|f\|_p^T.
$$
Since $\varepsilon$ is arbitrary then
$$
\int_C\widetilde{f}(y)dy\leq\|f\|_p^T,
$$
for each compact subset $C$ of $G$. Using the regularity Haar
measure we obtain $\widetilde{f}\in L^1(G)$. Since
$$
\int_G\widetilde{f}(y)dy=\int_G\Delta(y)^{-1/q}f(y)dy=\int_G\omega_q(y)f(y)dy,
$$
it follows that $f\in L^1(G,\omega_q)$ and so $f\in
L^1(G,\omega_q)\cap L^p(G)^+$. Also $L^1(G,\omega_q)\cap
L^p(G)^+\subseteq L^{T}_p(G)^+$, by Lemma \ref{l2}. Thus the proof
is completed.
\end{proof}
As a consequence of Theorem \ref{t4} we have the following
corollary. This result is indicated in \cite[Lemma 2.3]{f2} and
also follows from the main results in \cite{g1}.

\begin{cor}\label{c1}
Let $G$ be a unimodular amenable locally compact group and
$1<p\leq 2$. Then $L^1(G)\cap L^p(G)^+=L_p^{t}(G)^+=L_p^{T}(G)^+$.
\end{cor}

\begin{rem}\label{r1}\rm
Let $G$ be a locally compact group and $1<p<\infty$.
\begin{enumerate}
\item[(i)] Let $f\in L^{T}_p(G)$. Regarding to Remark \ref{r2}, it
is unknown to us whether the positive parts of $Re(f)$ and $Im(f)$
belong to $L^{T}_p(G)$. If this statement holds, for some locally
compact groups $G$ and also $1<p<\infty$, then Theorem \ref{t4}
implies that $L^1(G,\omega_q)\cap L^p(G)=L_p^{T}(G)$. Moreover
since $L^p(G)$ is a right $L^1(G,\omega_q)-$module, it follows
that $L^{T}_p(G)$ is an algebra under convolution. \item[(ii)] Let
$G$ be compact. We have already seen in Theorem \ref{t1} that
$L^p(G)=L^{T}_p(G)$. In fact
$$
L^1(G,\omega_q)\cap L^p(G)=L^1(G)\cap L^p(G)=L^p(G)=L^{T}_p(G).
$$
Consequently, part (i) of the current remark holds for compact
groups. \item [(iii)] Let $G$ be amenable. Then
$L^1(G)^+=L_p^{T}(G)^+$ if and only if $G$ is discrete. Indeed, if
$L^1(G)^+=L_p^{T}(G)^+$, then $L^1(G)\subseteq L^p(G)$ and so $G$
is discrete from \cite[Theorem 1]{z}. The converse is immediately
obtained by Corollary \ref{c1}. \item [(iv)] If $G$ is abelian,
Corollary \ref{c1} shows that $L^1(G)\cap L^p(G)^+=L_p^{T}(G)^+$.
Thus for a large class of nondiscrete and noncompact locally
compact groups, $L_p^{T}(G)^+$ is a proper subset of $L^1(G)^+$
and also $L^p(G)^+$. \item[(v)] Let $1\leq p_1\leq p_2<\infty$. By
\cite[Theorem 1]{z}, $L^{p_1}(G)^+=L^{p_2}(G)^+$ if and only if
$p_1=p_2$ or $G$ is finite. Whereas by part (iii), if $G$ is
discrete and amenable (not necessarily finite) then
$\ell^1(G)^+=L^{T}_p(G)^+$, for all $1\leq p<\infty$. In fact in
this case, all $L^{T}_p(G)^+$ coincide.
\end{enumerate}
\end{rem}

As an application of Remark \ref{r1}, we end this section with the
following examples.

\begin{exa}\rm
Let $1<p<\infty$.
\begin{enumerate}
\item [(1)] Let $\mathbb{T}$ be the multiplicative circle group.
Since $\mathbb{T}$ is a compact group, Remark \ref{r1} part (ii)
provides that $L^{T}_p(\mathbb{T})=L^p(\mathbb{T})$. \item[(2)]
Let $\mathbb{Z}$ be the additive group of integers. Then Remark
\ref{r1} part (iii) implies that
$L_p^{T}(\mathbb{Z})^+=\ell^1(\mathbb{Z})^+$. \item [(3)] Let
$\mathbb{R}$ be the additive group of real numbers. Since
$\mathbb{R}$ is abelian, then
$L_p^{T}(\mathbb{R})^+=L^1(\mathbb{R})\cap L^p(\mathbb{R})^+$, by
Remark \ref{r1} part (iv).

\end{enumerate}
\end{exa}

\section{\bf Applications of Plancherel Theorem on $L^{T}_2(G)$}

For an abelian locally compact group $G$, the dual group of $G$
will be denoted by $\widehat{G}$. A straightforward application of
the Plancherel Theorem is that if a $f\in L^1(G)\cap L^2(G)$, then
$\widehat{f}$, the Fourier transform of $f$, is in
$L^2(\widehat{G})$, and the Fourier transform map is an isometry
with respect to the $L^2-$norm. In other word, the Fourier
transform map restricted to $L^1(G)\cap L^2(G)$ has a unique
extension to a linear isometry from $L^2(G)$ onto
$L^2(\widehat{G})$. This unique extension also denoted by
$f\mapsto\widehat{f}$ and is called the Plancherel transformation.
Also $\widehat{f}$ is called the Plancherel transform of $f$.
Furthermore, the inverse Fourier transform
$L^2(\widehat{G})\rightarrow L^2(G)$, defined by
$f\mapsto\breve{f}$, is also a linear isometry and both
transformations are inverse of each other; see \cite[31.17]{hr}
and also \cite[31.18]{hr} for more information.

In this section we show that the restricted mapping of the
Plancherel Theorem to $L^{T}_2(G)\cap L^{\infty}(G)$, is also a
linear isometry onto $L^{T}_2(\widehat{G})\cap
L^{\infty}(\widehat{G})$, under the norm given by
$\|.\|_2^T+\|.\|_{\infty}$. It requires some preparations. It
again should be noted that in this situation, $L^t_2(G)=L^T_2(G)$.
Moreover it is an algebra under convolution product and
$\|.\|_2^t=\|.\|_2^T$; see \cite[Theorems 1.1 and 2.5]{f2}.

\begin{lemma}\label{l4}
Let $G$ be an abelian locally compact group. Then for all $f\in
L^{T}_2(G)$ and $g\in L^2(G)$,
$(f*g)^{\widehat{}}=\widehat{f}\widehat{g}$
\end{lemma}

\begin{proof}
Let $f\in L^{T}_2(G)$ and $g\in C_{00}(G)$. There exists the
sequence $(f_n)$ in $C_{00}(G)$ such that
$\lim_{n\rightarrow\infty}\|f-f_n\|_2=0$. Thus
$\lim_{n\rightarrow\infty}\|g*f-g*f_n\|_2=0$ and so the Plancherel
Theorem implies that
$\lim_{n\rightarrow\infty}\|(g*f)^{\widehat{}}-(g*f_n)^{\widehat{}}\|_2=0$.
Since $(g*f_n)^{\widehat{}}=\widehat{g}\widehat{f_n}$ and also
$$
\lim_{n\rightarrow\infty}\|\widehat{g}\widehat{f_n}-\widehat{g}\widehat{f}\|_2
\leq\lim_{n\rightarrow\infty}\|\widehat{g}\|_{\infty}\|f_n-f\|_2=0,
$$
then $(f*g)^{\widehat{}}=\widehat{f}\widehat{g}$. Now let $g\in
L^2(G)$ and $(g_n)$ be a sequence in $C_{00}(G)$ such that
$\lim_{n\rightarrow\infty}\|g-g_n\|_2=0$. Hence by the Plancherel
theorem we have
$$
\lim_{n\rightarrow\infty}\|(g*f)^{\widehat{}}-(g_n*f)^{\widehat{}}\|_2\leq
\lim_{n\rightarrow\infty}\|g_n-g\|_2\|f\|_2^T=0.
$$
Thus there exists a subsequence $(g_{n_k})$ of $(g_n)$ such that
the sequences $(\widehat{g}_{n_k})$ and $(g_{n_k}*f)^{\widehat{}}$
converge pointwise almost everywhere on $\widehat{G}$ to
$\widehat{g}$ and $(g*f)^{\widehat{}}$, respectively. These
observations show that for almost all $\chi\in\widehat{G}$,
\begin{eqnarray*}
(g*f)^{\widehat{}}(\chi)&=&\lim_{k\rightarrow\infty}(g_{n_k}*f)^{\widehat{}}(\chi)\\
&=&\lim_{k\rightarrow\infty}\widehat{g}_{n_k}(\chi)\widehat{f}(\chi)\\
&=&\widehat{g}(\chi)\widehat{f}(\chi).
\end{eqnarray*}
Consequently $(g*f)^{\widehat{}}=\widehat{g}\widehat{f}$.
\end{proof}

Parseval's identity \cite[31.19]{hr} states that for all $f,g\in
L^2(G)$,
$$
\int_Gf(x)\overline{g(x)}dx=\int_{\widehat{G}}\widehat{f}(\gamma)\overline{\widehat{g}(\gamma)}d\gamma.
$$
Also a generalization of this result establishes that
\begin{equation}\label{e8}
\langle f,\check{g}\rangle=\langle\widehat{f},g\rangle,
\end{equation}
where $f\in L^2(G)$ and $g\in L^2(\widehat{G})$. Another result of
the Parseval's identity is the following equation
\begin{equation}\label{e12}
(fg)^{\widehat{}}=\widehat{f}*\widehat{g},
\end{equation}
where $f,g\in L^2(G)$ \cite[31.29]{hr}. These observations
corroborate the next result.

\begin{lemma}\label{l3}
Let $G$ be an abelian locally compact group and $f,g\in
L^2(\widehat{G})$. Then $(fg)^{\check{}}=\check{f}*\check{g}$.
\end{lemma}
Let $f\in L^{\infty}(G)$. Then the map
$$
M_f:L^2(G)\rightarrow L^2(G)
$$
defined as $g\mapsto fg$, is well defined by \cite[Theorem
20.13]{hs} and belongs to $B(L^2(G))$, consisting of all bonded
linear operators from $L^2(G)$ into $L^2(G)$. Also with some
elementary calculation we have $\|M_f\|=\|f\|_{\infty}$. Note that
$L^{\infty}(G)$ can be considered as a Banach algebra under
pointwise multiplication of functions. It follows that the map
$f\mapsto M_f$ from $L^{\infty}(G)$ into $B(L^2(G))$ is an
isometric isomorphism. This fact will be used several times in the
current debate.
\begin{pro}\label{p6}
Let $G$ be an abelian locally compact group and $f\in L^2(G)$ and
$g\in L^2(\widehat{G})$.
\begin{enumerate}
\item[(i)] $f\in L^{T}_2(G)$ if and only if $\widehat{f}\in
L^{\infty}(\widehat{G})$ and $\|f\|_2^T=\|\widehat{f}\|_{\infty}$.
\item[(ii)] $g\in L^{T}_2(\widehat{G})$ if and only if
$\check{g}\in L^{\infty}(G)$ and
$\|g\|_2^T=\|\check{g}\|_{\infty}$
\end{enumerate}
\end{pro}

\begin{proof}
$(i)$. Let $\widehat{f}\in L^{\infty}(\widehat{G})$ and $g\in
C_{00}(G)$. Lemmas \ref{l3} and \ref{l4} together with the
Plancherel Theorem imply that
\begin{eqnarray*}
\|g*f\|_2&=&\|(g*f)^{\widehat{}}\|_2\\
&=&\|\widehat{g}\widehat{f}\|_2\\
&\leq&\|\widehat{g}\|_2\|\widehat{f}\|_{\infty}\\
&=&\|g\|_2\|\widehat{f}\|_{\infty}.
\end{eqnarray*}
It follows that $f\in L^{T}_2(G)$ and also
$\|f\|_2^T\leq\|\widehat{f}\|_{\infty}$. For the converse, let
$f\in L^{T}_2(G)$. Thus for each $g\in L^2(G)$,
$(f*g)^{\widehat{}}=\widehat{f}\widehat{g}$ by Lemma \ref{l4}. It
follows that $\widehat{f}\widehat{g}\in L^2(\widehat{G})$ and so
the linear operator $M_{\widehat{f}}:L^{2}(\widehat{G})\rightarrow
L^{2}(\widehat{G})$, defined by
$\widehat{g}\mapsto\widehat{f}\widehat{g}$ is well defined.
Moreover
\begin{eqnarray*}
\|M_{\widehat{f}}\|&=&\sup_{g\in
L^2(\widehat{G}),\|\widehat{g}\|_2\leq
1}\|\widehat{f}\widehat{g}\|_2\\
&=&\sup_{g\in L^2(G),\|g\|_2\leq 1}\|(g*f)^{\widehat{}}\|_2\\
&=&\sup_{g\in L^2(G),\|g\|_2\leq 1}\|g*f\|_2\\
&=&\|f\|_2^T.
\end{eqnarray*}
Therefore $\widehat{f}\in L^{\infty}(\widehat{G})$ and
$\|\widehat{f}\|_{\infty}=\|M_{\widehat{f}}\|$. One can prove the
statement $(ii)$ with some similar arguments.
\end{proof}

Now we are in a position to prove the main result of this section.

\begin{thm}\label{t5}
Let $G$ be an abelian locally compact group. Then the restricted
mapping of the Plancherel Theorem is a linear isometry from
$L^{T}_2(G)\cap L^{\infty}(G)$ onto $L^{T}_2(\widehat{G})\cap
L^{\infty}(\widehat{G})$. The inverse transformation is also a
linear isometry. Moreover These two transformations are inverse of
each other.
\end{thm}

\begin{proof}
Let $f\in L^{T}_2(G)\cap L^{\infty}(G)$. Then the Plancherel
Theorem and also Proposition \ref{p6} imply that $\widehat{f}\in
L^2(\widehat{G})\cap L^{\infty}(\widehat{G})$. That
$\widehat{f}\in L^{T}_2(\widehat{G})$ follows by Proposition
\ref{p6} part $(ii)$. Moreover, by (\ref{e12}) and Plancherel
Theorem we have
\begin{eqnarray*}
\sup_{h\in L^2(\widehat{G}),\|h\|_2\leq
1}\|h*\widehat{f}\|_2&=&\sup_{h\in L^2(\widehat{G}),\|h\|_2\leq
1}\|(\check{h}f)^{\widehat{}}\|_2\\
&=&\sup_{h\in L^2(\widehat{G}),\|h\|_2\leq
1}\|\check{h}f\|_2\\
&=&\sup_{g\in L^2(G),\|g\|_2\leq
1}\|gf\|_2\\
&=&\|f\|_{\infty}.
\end{eqnarray*}
It follows that $\|\widehat{f}\|^T_2=\|f\|_{\infty}$. Furthermore
it is provided from Proposition \ref{p6} that
$\|f\|^T_2=\|\widehat{f}\|_{\infty}$. Thus
$$
\|f\|^T_2+\|f\|_{\infty}=\|\widehat{f}\|_{\infty}+\|\widehat{f}\|^T_2
$$
and so the restricted mapping is isometric on $L^{T}_2(G)\cap
L^{\infty}(G)$. For the inverse transformation, suppose that $g\in
L^{T}_2(\widehat{G})\cap L^{\infty}(\widehat{G})$. Thus
$\check{g}\in L^2(G)$. With some similar arguments one can easily
show that $\check{g}\in L^{\infty}(G)$ and so $\check{g}\in
L^{T}_2(G)\cap L^{\infty}(G)$. Also
$$
\|\check{g}\|^T_2+\|\check{g}\|_{\infty}=\|g\|^T_2+\|g\|_{\infty}.
$$
Thus the restricted mappings are isometric. That two
transformations are inverse of each other follows from the
Plancherel Theorem.
\end{proof}

\vspace{9mm}

{\footnotesize \noindent
 F. Abtahi\\
  Department of Pure Mathematics,\\
  Faculty of Mathematics and Statistics,\\
   University of Isfahan,\\
    Isfahan 81746-73441, Iran\\
     f.abtahi@sci.ui.ac.ir,abtahif2002@yahoo.com\\

\noindent
 H. G. Amini\\
  Department of Mathematics,
   University of Isfahan,
    Isfahan, Iran\\
     heidaramini37@gmail.com\\

\noindent
 A. Rejali\\
  Department of Pure Mathematics,\\
  Faculty of Mathematics and Statistics,\\
   University of Isfahan,\\
    Isfahan 81746-73441, Iran\\
     rejali@sci.ui.ac.ir\\

\end{document}